\input amstex
\documentstyle{amsppt}
\magnification=\magstephalf
 \addto\tenpoint{\baselineskip 15pt
  \abovedisplayskip18pt plus4.5pt minus9pt
  \belowdisplayskip\abovedisplayskip
  \abovedisplayshortskip0pt plus4.5pt
  \belowdisplayshortskip10.5pt plus4.5pt minus6pt}\tenpoint
\pagewidth{6.5truein} \pageheight{8.9truein}
\subheadskip\bigskipamount
\belowheadskip\bigskipamount
\aboveheadskip=3\bigskipamount
\catcode`\@=11
\def\output@{\shipout\vbox{%
 \ifrunheads@ \makeheadline \pagebody
       \else \pagebody \fi \makefootline 
 }%
 \advancepageno \ifnum\outputpenalty>-\@MM\else\dosupereject\fi}
\outer\def\subhead#1\endsubhead{\par\penaltyandskip@{-100}\subheadskip
  \noindent{\subheadfont@\ignorespaces#1\unskip\endgraf}\removelastskip
  \nobreak\medskip\noindent}
\outer\def\enddocument{\par
  \add@missing\endRefs
  \add@missing\endroster \add@missing\endproclaim
  \add@missing\enddefinition
  \add@missing\enddemo \add@missing\endremark \add@missing\endexample
 \ifmonograph@ 
 \else
 \vfill
 \nobreak
 \thetranslator@
 \count@\z@ \loop\ifnum\count@<\addresscount@\advance\count@\@ne
 \csname address\number\count@\endcsname
 \csname email\number\count@\endcsname
 \repeat
\fi
 \supereject\end}
\catcode`\@=\active
\CenteredTagsOnSplits
\NoBlackBoxes
\nologo
\def\today{\ifcase\month\or
 January\or February\or March\or April\or May\or June\or
 July\or August\or September\or October\or November\or December\fi
 \space\number\day, \number\year}
\define\({\left(}
\define\){\right)}

\define\CC{{\Bbb C}}

\define\RR{{\Bbb R}}
\define\SS{\Bbb S}

\define\ZZ{{\Bbb Z}}
\define\[{\left[}
\define\]{\right]}
\define\ch{\operatorname{ch}}
\define\chiup{\raise.5ex\hbox{$\chi$}}
\define\cir{S^1}

\define\exertag #1#2{#2\ #1}

\define\index{\operatorname{index}}

\define\mstrut{^{\vphantom{1*\prime y}}}
\define\protag#1 #2{#2\ #1}

\define\res#1{\negmedspace\bigm|_{#1}}
\define\temsquare{\raise3.5pt\hbox{\boxed{ }}}

\define\theprotag#1 #2{#2~#1}

\define\xca#1{\removelastskip\medskip\noindent{\smc%
#1\unskip.}\enspace\ignorespaces }

\define\zmod#1{\ZZ/#1\ZZ}

\NoRunningHeads 
\define\Coker{\operatorname{Coker}}
\define\Fix{\operatorname{Fix}}
\define\Hor{Ho\v rava}
\define\Ker{\operatorname{Ker}}
\define\K{\Cal{K}}
\define\Pin{\operatorname{Pin}}
\define\Trace{\operatorname{Trace}}
\define\bX{\partial X}
\define\codim{\operatorname{codim}}
\define\hal{\frac12}
\define\psidot{\dot\psi}
\define\th{\hat\tau }
\define\ttw{\tilde{\tau }}
\define\zo{[0,1]}

\refstyle{A}
\widestnumber\key{SSSSS}   
\document

	\topmatter
 \title\nofrills  Two Index Theorems in Odd Dimensions \endtitle
 \author Daniel S. Freed  \endauthor
 \thanks The author is supported by NSF grant DMS-9307446, Presidential
Young Investigators award DMS-9057144, and by the O'Donnell
Foundation.\endthanks
 \affil Department of Mathematics \\ University of Texas at Austin\endaffil 
 \address Department of Mathematics, University of Texas, Austin, TX
78712\endaddress 
 \email dafr\@math.utexas.edu \endemail
 \date June 26, 1997\enddate
	\endtopmatter

\document

A recent paper of \Hor\ and Witten~\cite{HW}---part of the current flurry of
activity in string theory---contains an anomaly computation for $\cir/\langle
\tau \rangle\times \RR^{10}$, where $\langle \tau \rangle$~is the cyclic
group of order two generated by a reflection.  It was well established
10~years ago (e.g.~\cite{AS1}, ~\cite{F1}) that anomalies measure
nontriviality in the determinant line bundle of a family of Dirac operators,
and so can be computed topologically from the Atiyah-Singer index theory.
The novelty in the \Hor-Witten computation is a nontrivial index in {\it
odd\/} dimensions of a type not seen in standard index theory.  We abstract
two general theorems which imply the \Hor-Witten result.  (Naturally, we
replace~$\RR^{10}$ by a compact manifold~$Y^{10}$).  Theorem~A is a Lefschetz
formula for an {\it orientation-reversing\/} isometric involution on an odd
dimensional manifold.  The Atiyah-Bott-Segal-Singer applications of Lefschetz
theory~\cite{AB1}, ~\cite{ASe}, ~\cite{AS2} all deal with {\it
orientation-preserving\/} isometries for which there is no nontrivial
Lefschetz formula in odd dimensions~\cite{AS2, Proposition~9.3}.
Alternatively, we can consider $\zo\times Y$ in place of $S^1/\langle \tau
\rangle\times Y$ and then the \Hor-Witten anomaly computation is a boundary
value problem with {\it local\/} boundary conditions.  Theorem~B generalizes
this situation and is closely related to the boundary value problem in the
original proof of the Atiyah-Singer index theorem~\cite{P}.
 
Our proofs use standard techniques, except for a small trick used to prove
Theorem~A.  For simplicity we discuss the standard complex Dirac operator;
the theorems are true for any Dirac operator.  Our language refers mostly to
a single operator, though the results hold for families of Dirac operators as
required by the anomaly problem.  In this regard we remark that Theorem~A
only holds modulo 2-torsion in the $K$-theory of the parameter space, whereas
Theorem~B holds exactly in $K$-theory.  For the anomaly problem this means
that Theorem~A may not be adequate to detect all {\it global\/} anomalies.
(In the general situation of Theorem~A, there is probably no fixed-point
formula for the exact index.)
 
I thank Edward Witten for bringing this anomaly computation to my attention
and for discussions.

 \head
 \S{1} A Lefschetz formula for orientation-reversing isometries
 \endhead
 \comment
 lasteqno 1@ 21
 \endcomment

Let $X$~be a compact odd dimensional spin manifold.  Suppose $\tau \:X\to X$
is an orientation-reversing isometric involution.  Assume there exists a
lift~$\ttw\:S_X\to S_X$ to the complex spinor bundle~$S_X$ on~$X$ such that 
  $$ D_X\ttw = -\ttw D_X, \tag{1.1} $$ 
where $D_X$ is the Dirac operator.  It follows from \theprotag{1.5} {Lemma}
below that $\ttw^2$~is locally constant, so dividing by a square root of that
locally constant function we may assume
  $$  \ttw^2=1. \tag{1.2} $$
Then the $\pm1$-eigenspaces of~$\ttw$ give a splitting of the spinor fields
  $$ S(X) \cong S^+(X)\oplus S^-(X), \tag{1.3} $$ 
and the Dirac operator interchanges~$S^+(X)$ and~$S^-(X)$.  Our problem is to
compute 
  $$ \index\bigl[ D_X\: S^+(X)\longrightarrow S^-(X)\bigr]. \tag{1.4} $$ 
 
The simplest example is~$X=\cir=\RR/\ZZ$ with $\tau $~the
reflection~$x\mapsto -x$.  The spinor fields may be identified with the
complex functions, the Dirac operator with~$i\frac{d}{dx}$, and the
splitting~\thetag{1.3} is the splitting into even and odd functions.  Here
the index is~1.  The \Hor-Witten example is the product with a fixed even
dimensional manifold~$Y$, in which case the index is $\index D_Y$.

The lift~$\ttw$, if it exists, is almost unique.

        \proclaim{\protag{1.5} {Lemma}}
 Suppose $X$~is an odd dimensional spin manifold, and $\theta \:S_X\to S_X$ a
bundle map such that $D_X\theta =\theta D_X$.  Then $\theta $~is a locally
constant multiple of the identity.
        \endproclaim

\flushpar
 If $\ttw_1,\ttw_2$ are two lifts of~$\tau $ satisfying~\thetag{1.1}
and~\thetag{1.2}, set $\theta =\ttw_1\ttw_2$ to conclude that
$\ttw_1=\pm\ttw_2$ on each component of~$X$.

        \demo{Proof}
 Fix $x\in X$ and choose a local oriented orthonormal framing~$\{e_i\}$
near~$x$.  Then if $\psi $~is a spinor field with $\psi (x)=0$, an easy
computation shows
  $$ 0=(D_X\theta -\theta D_X)\psi (x) = \bigl[c(e^i),\theta (x)\bigr]\,\nabla
     _{e_i}\psi (x),  \tag{1.6} $$
where $\{e^i\}$~is the dual coframing, $c(\cdot )$~is Clifford
multiplication, and $\nabla $~is the Levi-Civita covariant derivative.  Fix
an index~$i$.  Choose a set of spinor fields~$\{\psi ^{(\alpha )}\}_{\alpha
}$ so that $\psi ^{(\alpha )}(x)=0$, the derivatives $\nabla _{e_j}\psi
^{(\alpha )}(x)=0 $ for~$j\not= i$, and $\{\nabla _{e_i}\psi ^{(\alpha)
}(x)\}_\alpha $ span the fiber~$(S_X)_x$.  Then~\thetag{1.6} implies
$\bigl[c(e^i),\theta (x)\bigr]=0$ for all~$i$, and since the spin
representation is irreducible in odd dimensions, $\theta (x)$~is a scalar.
Then for any spinor field~$\psi $,
  $$ 0=(D_X\theta -\theta D_X)\psi =c(d\theta )\psi ,  $$ 
from which $d\theta =0$ so that $\theta $~is locally constant. 
        \enddemo

Concerning the existence of~$\ttw$, we recall that in odd dimensions the spin
representation~$\SS$ extends to an ungraded module for the Clifford algebra
on which the volume form, suitably normalized, acts as~$+1$.  In particular,
$\SS$~is a representation of the Pin group.  Now the isometry~$\tau $ lifts
to the bundle of orthonormal frames~$O(X)$ of~$X$.  The spin structure
induces a pin structure~$\Pin(X)$---a principal Pin bundle which double
covers~$O(X)$---and it is a topological question about covering spaces to
determine if $\tau $~acting on~$O(X)$ lifts to~$\Pin(X)$.  If so, the lift
may have order~4.  In any case the spinor bundle~$S_X$ is associated
to~$\Pin(X)$, and the lift induces a map~$\ttw$ on spinor fields.  But
Clifford multiplication is {\it not\/} a map of Pin representations---there
is a sign for elements which reverse the orientation---and so the Dirac
operator does not extend simply extend to the Pin bundle.  Rather, the sign
means that the lift~$\ttw$ of an orientation-reversing isometry anticommutes
with the Dirac operator as in~\thetag{1.1}.
 
We turn now to the index~\thetag{1.4}.  The general Lefschetz formulas of
Atiyah-Bott-Segal-Singer~\cite{AB1}, \cite{ASe}, ~\cite{AS2} apply to an
elliptic operator $D\:C^{\infty}(E)\to C^{\infty}(F)$ acting between two
vector bundles~$E,F$ with endomorphisms $\th_E,\th_F$ such that 
  $$ D\th_E = \th_FD. \tag{1.7} $$ 
Our problem concerns the Dirac operator $D_X\:C^{\infty}(S_X)\to
C^{\infty}(S_X)$, but the given lift~$\ttw$ satisfies~\thetag{1.1},
not~\thetag{1.7}.  Here is the trick: Define
  $$ \th =\cases \hphantom{-}\ttw ,&\text{on the domain copy
     of~$S_X$};\\-\ttw,&\text{on the codomain copy of~$S_X$}.\endcases
     \tag{1.8} $$
Now $\th$ satisfies~\thetag{1.7}!  The Lefschetz number is 
  $$ \aligned
      L(\th,D_X) &= \Trace \th\res{\Ker D_X} - \Trace \th\res{\Coker D_X}, \\
      &= 2\,\index\bigl[ D_X\:S^+(X)\longrightarrow S^-(X)\bigr],\endaligned
      $$
twice the index we would like to compute.

The generalized Lefschetz formulas compute this index in terms of the fixed
point set~$\Fix(\tau )$ of~$\tau $.  In our situation each component~$F$
of~$\Fix(\tau )$ is an even dimensional manifold.  The Atiyah-Segal
formula~\cite{ASe,Theorem~2.12} applies in general; we first state the result
with the vastly simplifying assumption that the normal bundle~$N_F$ to each
component of~$\Fix(\tau )$ is trivial.  See \theprotag{1.10} {Remark}
following the statement of Theorem~A for the formula when $N_F$~is only
assumed orientable.

        \proclaim{Theorem~A}
 Let $X$~be an odd dimensional spin manifold, $\tau \:X\to X$ an
orientation-reversing isometric involution, and $\ttw\:S_X\to S_X$ a lift to
spinors which anticommutes with the Dirac operator~$D_X$ and
satisfies~$\ttw^2=1$.  Then $D_X$~exchanges the
$\pm1$-eigenspaces~$S^{\pm}(X)$ of~$\ttw$ operating on spinor fields.  Assume
that each component~$F$ of the fixed point set~$\Fix(\tau )$ has trivial
normal bundle.  The sum over these components appears in the index formula
  $$ \index\bigl[ D_X\: S^+(X)\longrightarrow S^-(X)\bigr] =
     \sum\limits_{F}\frac{\index D_F}{2^{r(F)+1}}. \tag{1.9} $$
Here $\codim F=2r(F)+1$ and $D_F\:S^+(F)\to S^-(F)$ is the chiral Dirac
operator on~$F$ relative to an orientation chosen compatibly with~$\ttw$.
        \endproclaim

\flushpar
 The orientation is explained in the proof (see~\thetag{1.20}).  We make
several remarks before proceeding to the proof.  

        \remark{\protag{1.10} {Remark}}
 More generally, suppose only that each component~$F$ of the fixed point set
has {\it orientable\/} normal bundle.  Then \thetag{1.9}~is replaced by 
  $$ \index\bigl[ D_X\: S^+(X)\longrightarrow S^-(X)\bigr] = \frac 12
     \sum\limits_{F} \frac{\hat\Cal{A}(F)}{\ch \Delta(N_F)}[F],
     \tag{1.11} $$
where $\hat\Cal{A}$~is the usual characteristic class associated to Dirac,
$\ch \Delta$~is the Chern character of the spin bundle, and the orientation
of~$F$ is determined below.  (One does not need a spin structure to
define~$\ch\Delta $.)  See~~\cite{AS2,\S5} for a similar result.  This
formula only holds rationally in families.
        \endremark

        \remark{\protag{1.12} {Remark}}
 For the \Hor-Witten example $X=\cir\times Y$, $\tau $~is reflection on the
$\cir$ factor, and Theorem~A computes
  $$ \index\bigl[ D_X\: S^+(X)\longrightarrow S^-(X)\bigr] = \index D_Y,
     \tag{1.13} $$
which agrees with~\cite{HW}.  Here $Y$~is a compact even dimensional spin
manifold.  According to \theprotag{1.15} {Remark} below this only holds
modulo 2-torsion in families.  In the next section we show that in fact this
result holds exactly (see~\thetag{2.11}).
        \endremark

        \remark{\protag{1.14} {Remark}}
 For a single operator we can use the heat kernel approach to the Lefschetz
formula (see~\cite{R}, ~\cite{BGV} for example) to derive~\thetag{1.9}.  We
write 
  $$ \index\bigl[ D_X\: S^+(X)\longrightarrow S^-(X)\bigr] = \int_{X\times
     X}\Trace\bigl(\ttw(x,y)e^{-tD_X^2}(y,x) \bigr)\;dy\,dx,  $$
valid for any~$t$, and let~$t\to 0$.  The integral then localizes on the
fixed point set.  As always in index theory, this heat kernel  approach does
not generalize to the integral $K$-theory index of a family of Dirac
operators.
        \endremark

        \remark{\protag{1.15} {Remark}}
 Theorem~A applies to families of Dirac operators, but only gives a result
in~$K(Z)[\frac 12]$, where $Z$~is the parameter space.  (Below we
use~\cite{ASe,Theorem~2.12}.  Although this theorem is stated for
$K\text{-theory}\otimes \CC,$ in our situation the localization of the global
symbol of Dirac only involves denominators which are powers of~2.)
        \endremark

        \remark{\protag{1.16} {Remark}}
 Theorem~A also applies to (families of) real Dirac operators and Dirac
operators coupled to other vector bundles.  The \Hor-Witten example is
actually for the {\it real\/} Dirac operator coupled to the tangent bundle.
The quantity of interest is the {\it square root\/} of the determinant line
bundle, which is computed in $KO$-theory.  (See~\cite{F2,\S3} for an
explanation of this square root.)
        \endremark

        \demo{Proof of Theorem~A}
 We apply~\cite{ASe,Theorem~2.12} which asserts 
  $$ L(\th,D_X) = \sum\limits_{F}\index\left\{ \frac{\iota _F^*\,\sigma
     (D_X)(\th)}{\lambda _{-1}(N_F\otimes \CC)(\th)} \right\} , \tag{1.17} $$
where $\iota _F\:F\hookrightarrow X$ is the inclusion, $\sigma (D_X)\in
K_G(TX)$ is the symbol of Dirac, and 
  $$ \lambda _{-1}(N_F\otimes \CC) = \sum\limits_{}(-1)^i{\tsize\bigwedge}
     ^i(N_F\otimes \CC)\in K_G(F).  $$
Here $G=\langle \th  \rangle$, the cyclic group generated by~$\th$.
Evaluation on~$\th$ is the homomorphism 
  $$ K_G(F)@>\cong >> K(F)\otimes R(G)@>>> K(F)  $$ 
which evaluates a virtual character on~$\th$.  (For the cyclic group of order
two the virtual characters are real-valued.)
     
We work on a fixed component~$F$ of codimension~$2r(F)+1=2r+1$.  Since
$N_F$~is assumed trivial, we have an isomorphism $N_F\cong L^{\oplus (2r+1)}$
in~$K_G(F)$, where $L$~is the trivial real line bundle with $\th$~acting
as~$-1$.  It follows easily that 
  $$ \lambda _{-1}(N_F\otimes \CC)(\th)=2^{2r+1}. \tag{1.18} $$ 
 
Recall that the symbol $\sigma (D_X)$ evaluated on a cotangent vector~$\theta
$ is Clifford multiplication $c(\theta )\: S_X\to S_X$.  We need to compute
this for $\theta $~a cotangent vector to~$F$.  First, note that $F$~is
orientable, since $N_F$~is trivial.  We fix the orientations of~$F$ and~$N_F$
below.  Let $N_F$~have the trivial spin structure.  This, together with the
spin structure on~$X$, induces a spin structure on~$F$.  Then, letting
$S_F,S_{N_F}$ denote the spin bundles on the tangent and normal bundles
to~$F$, we have
  $$ S_X\res F\cong S_F\otimes S_{N_F}\cong S_F^{\oplus (2^r)} \tag{1.19} $$ 
since the normal bundle is trivial.  Therefore, $\iota _F^*\sigma (D_X)$~is
Clifford multiplication on $2^r$~copies of~$S_F$.
 
To compute the action of~$\th$ we fix an equivariant tubular neighborhood
of~$F$, which is diffeomorphic to~$F\times \RR^{2r+1}$, and introduce a
product metric.  (This computation is local, so does not use the triviality
of~$N_F$.)  Let $e_i$ be the standard orthonormal basis of~$\RR^{2r+1}$,
$x^i$~the standard coordinates on~$\RR^{2r+1}$, and $f$ a coordinate on~$F$.
Then
  $$ D_X = D_F + c(e^i)\nabla _{e_i}.$$
We claim
  $$ \aligned \tau (f;x^1,\dots,x^{2r+1}) &= \langle f;-x^1,\dots,-x^{2r+1}
     \rangle \\
      (\ttw\psi )(f;x^1,\dots,x^{2r+1}) &= \pm i^{r+1}c(e^1)\dots
     c(e^{2r+1})\psi (f;-x^1,\dots,-x^{2r+1}) \\
      &= i^mc(\omega _F)\psi (f;-x^1,\dots,-x^{2r+1}), \endaligned
     \tag{1.20} $$
where $\psi $~is a spinor field, $\dim F=2m$, and $\omega _F$ is a real
volume form on~$F$ with $c(\omega _F)^2=(-1)^m$.  A routine computation shows
that the first expression for~$\ttw$ satisfies~\thetag{1.1} and~$\ttw^2=1$,
whence the remark following \theprotag{1.5} {Lemma} implies that this is the
correct expression (with one of the signs).  The second expression for~$\ttw$
follows from a simple computation with Clifford algebras.  It
determines~$\omega _F$ uniquely.
      
Now we fix the orientation on~$F$ so that $\omega _F$~is an oriented volume
form.  Then $S^\pm_F$~are the $\pm(i^{-m})$-eigenspaces of~$c(\omega _F)$
acting on~$S_F$, which by~\thetag{1.20} are the $\pm1$-eigenspaces of~$\ttw$.
Use \thetag{1.8}~and \thetag{1.19}~to conclude that
  $$ \iota _F^*\,\sigma (D_X)(\th) = 2^{r+1}\,\sigma (D_F). \tag{1.21} $$ 
The desired result~\thetag{1.9} follows from~\thetag{1.17}, \thetag{1.21},
and~ \thetag{1.18}. 
        \enddemo

 \newpage
 \head
 \S{2} An index theorem for manifolds with boundary
 \endhead
 \comment
 lasteqno 2@ 20
 \endcomment

Let~$X$ be a compact odd dimensional spin manifold with boundary.  The
orientation on~$X$ determines an orientation on~$\bX$ and so a splitting 
  $$ S_X\res{\bX}\cong S\mstrut _{\bX}\cong S^+_{\bX} \oplus S^-_{\bX}
     \tag{2.1} $$
of the spin bundle on the boundary.  This splitting leads to local boundary
conditions~$P^{\pm}$ for the Dirac operator~$D_X$: the domain
of~$(D_X,P^{\pm})$ is the set of spinor fields~$\psi $ on~$X$ with 
  $$  \bigl(\psi \res{\bX} \bigr)^{\pm}=0,  $$
where $\phi =\phi ^++\phi ^-$ is the
decomposition of a spinor field~$\phi \in S(\bX)$ relative to~\thetag{2.1}.
These local boundary value problems are a key ingredient in the original
proof of the Atiyah-Singer index theorem.  Indeed~\cite{P,\S17},
~\cite{BW,\S21}
  $$ \index(D_X,P^{\pm})=0. \tag{2.2} $$ 
This is used to show that the index of the chiral Dirac operator on the
boundary vanishes:
  $$ \index D_{\bX}=0. \tag{2.3} $$ 
Equation~\thetag{2.3} is the assertion that the index is a bordism
invariant. 
 
We consider a mixture of these boundary conditions.  Namely, we independently
choose~$P^+$ or~$P^-$ on each component of the boundary. 

        \proclaim{Theorem~B}
 Let $X$~be a compact odd dimensional spin manifold with boundary, and
$\bX=\sqcup_i Y_i$ the decomposition of the boundary into components.
For each~$i$ choose $\epsilon _i=+$ or $\epsilon _i=-$ and consider the Dirac
operator~$(D_X,P^\epsilon )$ whose domain is the set of spinor fields~$\psi $
such that 
  $$ \bigl(\psi \res{Y_i} \bigr)^{\epsilon _i} = 0. \tag{2.4} $$ 
Then 
  $$ \index(D_X,P^\epsilon ) = \sum\limits_{\ssize{i \text{
     with}}\atop{\ssize\epsilon _i=-}} \index D_{Y_i}= -\sum\limits_{\ssize{i
     \text{ with}}\atop{\ssize\epsilon _i=+}} \index D_{Y_i} \tag{2.5} $$
        \endproclaim

\flushpar
 Note that the last equality follows directly from~\thetag{2.3}.  Also, if
all~$\epsilon _i=+$ or all~$\epsilon _i=-$, then \thetag{2.5}~reduces
to~\thetag{2.2} in view of~\thetag{2.3}.  As is evident from the proof below,
Theorem~B is a direct consequence of well-known facts about boundary-value
problems for Dirac operators.

        \remark{\protag{2.6} {Remark}}
 Theorem~B also holds in families; then \thetag{2.5} is an {\it exact\/}
equation in~$K(Z)$, where $Z$~is the parameter space.  (Contrast with
Theorem~A which only holds in~$K(Z)[\frac 12]$.)  As with Theorem~A (see
\theprotag{1.16} {Remark}), Theorem~B holds for (families of) real Dirac
operators and Dirac operators coupled to other vector bundles.
        \endremark

        \remark{\protag{2.7} {Remark}}
 Consider $X=[0,\hal]\times Y$, where $Y$ is a closed even dimensional spin
manifold.  We use the product metric.  Then $\bX = Y_0\sqcup Y_{\hal}$, where
$Y_{\hal}\cong Y$ and $Y_0\cong -Y$.  Here `$-Y$' denotes $Y$~with the
opposite orientation.  Let $\epsilon _0=+$ and $\epsilon _{\hal}=-$.  Then
\thetag{2.5} gives 
  $$ \index(D_X,P^\epsilon ) = \index D_Y. \tag{2.8} $$ 
Let $\tilde{D}\:H^+\to H^-$ be the Dirac operator in the \Hor-Witten
example~\thetag{1.13}.  Here we are working on~$\cir\times Y$ and $\psi \in
H^{\pm}=S^{\pm}(\cir\times Y)$ is an $S(Y)$-valued function on~$\cir=\RR/\ZZ$
satisfying
  $$ \psi (-x) = \pm i^mc(\omega _Y)\psi (x), \tag{2.9} $$ 
where $\omega _Y$~is a volume form on~$Y$ and $\dim Y=2m$
(cf.~\thetag{1.20}).  We now give an {\it a priori\/} argument that 
  $$ \index(D_X,P^\epsilon ) = \index \tilde{D}, \tag{2.10} $$ 
even in families.  This is consistent with the computations~\thetag{1.13}
and~\thetag{2.8} from Theorem~A and Theorem~B for single operators, and gives
the exact result 
  $$ \index \tilde{D}=\index D_Y \tag{2.11} $$ 
in families.  (This was previously proved modulo 2-torsion.) 
 
To prove~\thetag{2.10} note that relative to the splitting $S(Y)\cong
S^+(Y)\oplus S^-(Y)$ equation~\thetag{2.9} asserts that $\psi \in
H^{\pm}$~satisfies
  $$ \aligned
      \psi ^+(-x) &= \pm \psi ^+(x), \\ 
      \psi ^-(-x) &= \mp \psi ^-(x) .\endaligned  $$ 
Consider the diagram 
  $$ \minCDarrowwidth{25pt}\CD
      0 @>>> C^{1,\delta }\bigl[H^+ \bigr] @>>> C^{1,\delta
     }\bigl[S([0,\hal]\times Y,\epsilon ) \bigr] @>p>> C^{\delta
     }\bigl[S^+(Y) \bigr] \oplus C^\delta\bigl[ S^+(Y)\bigr] @>>> 0\\
      @. @VV\tilde{D}V @VV{(D_X,\epsilon )}V @VV\operatorname{id}V\\
      0 @>>> C^\delta\bigl[H^- \bigr] @>>> C^{\delta
     }\bigl[S([0,\hal]\times Y) \bigr] @>q>> C^\delta\bigl[S^+(Y) \bigr]
     \oplus C^\delta\bigl[ S^+(Y)\bigr] @>>> 0 \endCD $$
where `$C^{1,\delta }[\cdot ]$' and `$C^\delta [\cdot ]$' denote spaces of
H\"older functions for some $0<\delta <1$; `$S([0,\hal]\times Y,\epsilon )$'
denotes the space of spinor fields satisfying~\thetag{2.4}, which in this
case is $\psi ^-(0) = \psi ^-(\hal)=0$; the first horizontal arrows are
restriction maps; and
  $$ \aligned
      p(\psi ) &= \langle -i\psidot^+(0), -i\psidot^+(\hal) \rangle, \\
      q(\psi ) &= \langle \psi^+(0), \psi^+(\hal) \rangle.\endaligned
      $$
(Here $\psidot=\frac{d\psi }{dx}$.)  A routine check shows that the rows are
exact and the diagram commutes.  Now \thetag{2.10}~is a consequence of the
following lemma.  (See~\cite{S} for a more general discussion.)

        \proclaim{\protag{2.12} {Lemma}}
 Let 
  $$ \CD 
      0 @>>> V' @>>> V @>>> V'' @>>> 0 \\
      @. @VVT_z'V @VVT_zV @VVT''_zV\\ 
      0 @>>> W' @>>> W @>>> W'' @>>> 0  
     \endCD \tag{2.13} $$ 
be a commutative diagram with exact rows, where $V',V,V'',W',W,W''$ are
Banach spaces and $T'_z,T\mstrut _z,T''_z$ are Fredholm operators depending
continuously on a parameter~$z\in Z$.  Then
  $$ \index(T) = \index(T' ) + \index(T'' )\in
     K(Z). \tag{2.14} $$
        \endproclaim

\flushpar
 The Banach spaces are allowed to vary continuously; we omit this from the
notation for convenience. 

        \demo{Proof}
 The short exact sequence of chain complexes~\thetag{2.13} induces a long
exact sequence in cohomology: 
  $$ 0 @>>> \Ker T'_z @>>> \Ker T_z @>>> \Ker T''_z @>>> \Coker T'_z @>>>
     \Coker T_z @>>> \Coker T''_z @>>> 0. \tag{2.15} $$
The exactness of~\thetag{2.15} proves~\thetag{2.14} for a single operator.
For a family it suffices to prove~\thetag{2.14} for $Z$~compact.
Then~\cite{AS3,\S2} we can find $w_1'(z),\dots ,w'_{N'}(z)\in W'$ and
$w_1(z),\dots ,w_N(z)\in W$ so that 
  $$ \CD 
      0 @>>> V' \oplus \CC^{N'}@>>> V \oplus \CC^{N'} \oplus  \CC^{N}@>>> V''
     \oplus \CC^N @>>> 0 \\
      @. @VVS_z'V @VVS_zV @VVS''_zV\\ 
      0 @>>> W' @>>> W @>>> W'' @>>> 0  
     \endCD   $$ 
satisfies the hypotheses of the lemma and in addition $S_z', S\mstrut
_z,S_z''$~are surjective.  Here
  $$ S_z(v;\lambda ^i;\mu ^j) = T_z(v) + \lambda ^iw_i'(z) + \mu ^jw_j(z)
      $$
and $S'_z,S''_z$~are the corresponding induced maps.  We have 
  $$ \index S' = \index T' + [Z\times \CC^{N'}] \in
     K(Z) $$
with similar formulas for the indices of~$S$ and~$S''$.  Now the exactness
of~\thetag{2.15} (with all cokernels vanishing) proves~\thetag{2.14}; the
extra trivial bundles cancel out.
        \enddemo
        \endremark

        \demo{Proof of Theorem~B}
 The proof is based on analysis by Calder\'on and Seeley~\cite{P,\S17}; we
rely on the account in~\cite{BW}.  We remark that the index with {\it
local\/} boundary conditions is a topological invariant; in fact, it has an
interpretation in $K$-theory~\cite{AB2}.  So, for example, we can deform the
metric to a metric which is a product near the boundary.
 
Consider first a single operator.  Let 
  $$ \hat{\K} = \Ker \bigl[ D_X\:S(X)\longrightarrow S(X) \bigr] $$ 
and $\K\subset S(\bX)$ the image of~$\hat{\K}$ under restriction to the
boundary.  We use the Sobolev completions~$H^1$ of~$S(X)$ and $H^{1/2}$
of~$S(\bX)$.  Then $\K$~is a closed infinite dimensional subspace
of~$S(\bX)$.  Let 
  $$ P^{\epsilon }\: S(\bX)\longrightarrow \bigoplus _i S^{\epsilon _i}(Y_i)
      $$
be the projection defined by the boundary condition~\thetag{2.4}.  The first
result~\cite{BW,Theorem~20.12} is that 
  $$ \index(D_X,P^\epsilon ) = \index\bigl[ P^\epsilon_{\K}
     \:\K\longrightarrow \bigoplus_i S^{\epsilon _i}(Y_i)
     \bigr], \tag{2.16} $$
where `$P^\epsilon _{\K}$'~denotes the restriction of~$P^\epsilon $ to~$\K$.
This applies in particular to~$P^+$ (which is~$P^\epsilon $ with all
$\epsilon _i=+$), and so~\cite{BW,Theorem~21.2}
  $$ \split
      \index(D_X,P^\epsilon ) - \index(D_X,P^+) &= \index(P^\epsilon _{\K}) -
     \index(P^+_{\K}) \\
      &= \index(P^\epsilon _{\K}) + \index\bigl(P^+_{\K})^*  \\
      &= \index\bigl[ P^\epsilon _{\K}(P^+_{\K})^*\: \bigoplus\limits_
     {\ssize{i \text{ with}}\atop{\ssize\epsilon _i=-}}
     S^+(Y_i)\longrightarrow \bigoplus\limits_{\ssize{i \text{
     with}}\atop{\ssize\epsilon _i=-}} S^-(Y_i) \bigr].\endsplit \tag{2.17} $$
The final step is the assertion (see~\cite{BW,Theorem~21.5}) that $P^\epsilon
_{\K}(P^+_{\K})^*$~is a pseudodifferential operator of order~0 whose
symbol---up to a factor and after restriction to the sphere bundle---is the
symbol of the Dirac operator $\sum\limits_{{\ssize{i \text{
with}}\atop{\ssize\epsilon _i=-}}}\!\!\!D_{Y_i}$.  (This is the brunt of the
argument; it depends on properties of the {\it Calder\'on projector\/}.)
Then the first equality in~\thetag{2.5} follows directly from~\thetag{2.17}
and~\thetag{2.2}.
      
We briefly consider how to modify this argument for a family of Dirac
operators parameterized by~$z\in Z$.  It suffices to consider $Z$~compact for
index computations.  Then as in the proof of the lemma above we can find a
finite number of spinor fields $\psi _1(z),\dots ,\psi _N(z)$ so that
  $$ \aligned
      \bigl(T(z),P^\epsilon (z) \bigr)\: S_{P^\epsilon (z)}(X)\oplus \CC^N
     &\longrightarrow S(X) \\
      \langle \psi ;\lambda ^i \rangle&\longmapsto
     D_X(z)\psi + \lambda ^i\psi _i(z)\endaligned \tag{2.18} $$
is surjective.  Here $S_{P^\epsilon (z)}(X)\subset S(X)$ is the subspace of
spinor fields satisfying the boundary condition~$P^\epsilon (z)$.  Then
  $$ \index(T,P^\epsilon) = \index(D_X,P^\epsilon ) + [Z\times \CC^N]\in
     K(Z). \tag{2.19} $$
Now $T(z)$---the operator~\thetag{2.18} extended to all of $S(X)\oplus
\CC^N$---is also surjective.  Thus the orthogonal complement to the kernel
of~$T(z)$ varies continuously in~$z$ (using $T(z)$~as an isomorphism to the
continuously varying codomains), whence the kernel $\hat\K(z)$ of~$T(z)$
varies continuously as well.  So does its image~$\K(z)$ in~$S(\bX)$
Equation~\thetag{2.16} is replaced by
  $$ \index(T,P^\epsilon ) = \index\bigl[P^\epsilon \:\K\longrightarrow
     \bigoplus_i S^{\epsilon _i}(Y_i) \bigr]. \tag{2.20} $$
This follows simply by identifying the kernel bundle of the families of
operators on each side; the cokernels vanish.  By adding more~$\psi _i(z)$ we
can ensure that \thetag{2.18}~is also surjective for~$\bigl(T(z),P^+(z)
\bigr)$ and repeat~\thetag{2.19} and~\thetag{2.20} for~$P^+$
replacing~$P^\epsilon $.  Then equation~\thetag{2.17} holds---the auxiliary
trivial bundle cancels out---and the proof concludes as before.
        \enddemo

\newpage
\widestnumber\key{SSSSSSS}   

\enddocument